\documentclass[10pt,letterpaper]{article}
\usepackage[top=0.85in,left=2.75in,footskip=0.75in]{geometry}

\usepackage{changepage}

\usepackage[utf8x]{inputenc}

\usepackage{textcomp,marvosym}

\usepackage{fixltx2e}

\usepackage{amsmath,amssymb}

\usepackage{cite}

\usepackage[pdftex]{graphicx}    
\usepackage{subfigure}   

\usepackage[right]{lineno}
\usepackage{subfigure}
\usepackage{amssymb,amsmath}
\usepackage{relsize}
\usepackage{multirow}

\usepackage{microtype}
\DisableLigatures[f]{encoding = *, family = * }

\raggedright
\usepackage[aboveskip=1pt,labelfont=bf,labelsep=period,justification=raggedright,singlelinecheck=off]{caption}

\makeatletter
\renewcommand{\@biblabel}[1]{\quad#1.}
\makeatother

\date{}

\usepackage{lastpage,fancyhdr,graphicx}
\usepackage{epstopdf}
\fancyheadoffset[L]{2.25in}
\fancyfootoffset[L]{2.25in}
\begin{document}
\vspace*{0.2in}

\begin{flushleft}
{\Large
\textbf\newline{Hybrid Modeling and Prediction of Dynamical Systems} 
}
\newline
\\
Franz Hamilton\textsuperscript{1,2*},
Alun Lloyd \textsuperscript{1,2,3},
Kevin Flores \textsuperscript{1,2,4}
\\
\bigskip
\textbf{1} Department of Mathematics, North Carolina State University, Raleigh, NC, USA
\\
\textbf{2} Center for Quantitative Sciences in Biomedicine, North Carolina State University, Raleigh, NC, USA
\\
\textbf{3} Biomathematics Graduate Program, North Carolina State University, Raleigh, NC, USA
\\
\textbf{4} Center for Research in Scientific Computation, North Carolina State University, Raleigh, NC, USA
\bigskip

%
%

* fwhamilt@ncsu.edu

\end{flushleft}
\section*{Abstract}
Scientific analysis often relies on the ability to make accurate predictions of a system's dynamics. Mechanistic models, parameterized by a number of unknown parameters, are often used for this purpose. Accurate estimation of the model state and parameters prior to prediction is necessary, but may be complicated by issues such as noisy data and uncertainty in parameters and initial conditions. At the other end of the spectrum exist nonparametric methods, which rely solely on data to build their predictions. While these nonparametric methods do not require a model of the system, their performance is strongly influenced by the amount and noisiness of the data. In this article, we consider a hybrid approach to modeling and prediction which merges recent advancements in nonparametric analysis with standard parametric methods. The general idea is to replace a subset of a mechanistic model's equations with their corresponding nonparametric representations, resulting in a hybrid modeling and prediction scheme. Overall, we find that this hybrid approach allows for more robust parameter estimation and improved short-term prediction in situations where there is a large uncertainty in model parameters. We demonstrate these advantages in the classical Lorenz-63 chaotic system and in networks of Hindmarsh-Rose neurons before application to experimentally collected structured population data.

\section*{Author Summary}
The question of how best to predict the evolution of a dynamical system has received substantial interest in the scientific community. While traditional mechanistic modeling approaches have dominated, data-driven approaches which rely on data to build predictive models have gained increasing popularity. The reality is, both approaches have their drawbacks and limitations. In this article we ask the question of whether or not a hybrid approach to prediction, which combines characteristics of both mechanistic modeling and data-driven modeling, can offer improvements over the standalone methodologies. We analyze the performance of these methods in two model systems and then evaluate them on experimentally collected population data.


\section*{Introduction}
Parametric modeling involves defining an underlying set of mechanistic equations which describe a system's dynamics. These mechanistic models often contain a number of unknown parameters as well as an uncertain state, both of which need to be quantified prior to use of the model for prediction. The success of parametric prediction is tied closely to the ability to construct accurate estimates of the model parameters and state. This can be particularly challenging in high dimensional estimation problems as well as in chaotic systems \cite{voss,baake}. Additionally, there is often a degree of model error, or a discrepancy between the structure of the model and that of the system, further complicating the estimation process and hindering prediction accuracy.

Despite these potential issues, mechanistic models are frequently utilized in data analysis. The question we aim to address is when is it advantageous to use them? Under suitable conditions where model error is relatively small and parameters can be reliably estimated, parametric predictions can provide a great deal of accuracy. However, as we will see in the subsequent examples, a large uncertainty in the initial parameter values often leads to inaccurate estimates resulting in poor model-based predictions.

An alternative approach to modeling and prediction abandons the use of any mechanistic equations, instead relying on predictive models built from data. These nonparametric methods have received considerable attention, in particular those methods based on Takens' delay-coordinate method for attractor reconstruction \cite{farmer,casdagli1989nonlinear,Sugihara:1990aa,smith1992identification,jimenez1992forecasting,sauer94,sugihara1994nonlinear,schroer1998predicting,kugiumtzis1998regularized,yuan,hsieh2005distinguishing,strelioff2006medium,regonda,schelter2006handbook,hamilton2016}.  The success of nonparametric methods is strongly influenced by the amount of data available as well as the dimension of the dynamical system. If only a sparse amount of training data is available, the result is often inaccurate predictions due to the lack of suitable nearby neighbors in delay-coordinate space. Furthermore, as the dimension and complexity of the dynamical system increases, nonparametric prediction becomes significantly more difficult due to the necessary data requirements \cite{hamilton2016}.

Several recent works have investigated the situation where only a portion of a mechanistic model is known \cite{hamilton2,berry2016}. Our motivation here though is to explore how best to use a full mechanistic model when it is available. We consider a hybrid methodology to modeling and prediction that combines the complementary features of both parametric and nonparametric methods. In our proposed hybrid method, a subset of a mechanistic model's equations are replaced by nonparametric evolution. These nonparametrically advanced variables are then incorporated into the remaining mechanistic equations during the data fitting and prediction process. The result of this approach is a more robust estimation of model parameters as well as an improvement in short-term prediction accuracy when initial parameter uncertainty is large.

The utility of this method is demonstrated in several example systems. The assumption throughout is that noisy training data from a system are available as well as a mechanistic model that describes the underlying dynamics. However, several of the model parameters are unknown and the model state is uncertain due to the noisy measurements. The goal is to make accurate predictions of the system state up to some forecast horizon beyond the end of the training data. We compare the prediction accuracy of the standard parametric and nonparametric methodologies with the novel hybrid method presented here.

We begin our analysis by examining prediction in the classical Lorenz-63 system \cite{lorenz63}, which exhibits chaotic dynamics. Motivated by the success of the hybrid method in the Lorenz-63 system, we consider a more sophisticated example of predicting the spiking dynamics of a neuron in a network of Hindmarsh-Rose \cite{hindmarsh} cells. Finally, we examine the prediction problem in a well-known experimental dataset from beetle population dynamics \cite{constantino}.

\section*{Materials and Methods}
The assumption throughout is that a set of noisy data is available over the time interval $\left[t(0),t(T)\right]$. This is referred to as the {\it training data} of the system. Using these training data, the question is how best to predict the system dynamics over the interval $\left[t(T+1),t(T+T_F)\right]$, known as the {\it prediction interval}. Standard parametric and nonparametric methods are presented before our discussion of the novel hybrid method which blends the two approaches. 

\subsection*{Parametric Modeling and Prediction}
When a full set of mechanistic equations is used for modeling and prediction, we refer to this as the parametric approach. Assume a general nonlinear system of the form
\begin{eqnarray} \label{e1}
\mathbf{x}(k+1) &=& \mathbf{f}\left(t(k),\mathbf{x}(k),\mathbf{p}\right)+\mathbf{w}(k)\\
\mathbf{y}(k) &=& \mathbf{h}\left(t(k),\mathbf{x}(k),\mathbf{p}\right)+\mathbf{v}(k)\nonumber
\end{eqnarray}
where $\textbf{x}= \left[x_1,x_2,\hdots,x_n\right]^{\mathsmaller T}$ is an $n$-dimensional vector of model state variables and $\textbf{p} = \left[p_1,p_2,\hdots,p_l\right]^{\mathsmaller T}$ is an $l$-dimensional vector of model parameters which may be known from first principles, partially known or completely unknown. $\textbf{f}$ represents our system dynamics which describe the evolution of the state $\mathbf{x}$ over time and \textbf{h} is an observation function which maps \textbf{x} to an $m$-dimensional vector of model observations, $\textbf{y} = \left[y_1,y_2,\hdots,y_m\right]^{\mathsmaller T}$. To simplify the description of our analysis, we assume that the training data maps directly to some subset of $\mathbf{x}$. $\textbf{w}(k)$ and $\textbf{v}(k)$ are assumed to be mean $\mathbf{0}$ Gaussian noise terms with covariances $\mathbf{Q}$ and $\mathbf{R}$ respectively. While discrete notation is used in Eq. \ref{e1} for notational convenience,  the evolution of \textbf{x} is often described by continuous-time systems. In this situation numerical solvers, such as Runge-Kutta or Adams-Moulton methods, are used to obtain solutions to the continuous-time system at discrete time points.

When the state of a system is uncertain due to noisy or incomplete observations, nonlinear Kalman filtering can be used for state estimation \cite{voss,enkf7,evensen,rabier,cummings,yoshida,stuart,schiffbook,berry2,hamiltonEPL,hamiltonPRE,ghanim,ghanim2,sitz2002}. Here we choose the unscented Kalman filter (UKF), which approximates the propagation of the mean and covariance of a random variable through a nonlinear function using a deterministic ensemble selected through the unscented transformation \cite{simon,julier1,julier2}. We initialize the filter with state vector $\mathbf{x^{+}}(0)$ and covariance matrix $\mathbf{P^{+}}(0)$. At the $k$th step of the filter there is an estimate of the state $\mathbf{x^{+}}(k-1)$ and the covariance matrix $\mathbf{P^{+}}(k-1)$. In the UKF, the singular value decomposition is used to find the square root of the matrix $\mathbf{P^{+}}(k-1)$, which is used to form an ensemble of $2n+1$ state vectors.

The model $\mathbf{f}$ is applied to the ensemble, advancing it forward one time step, and then observed with $\mathbf{h}$. The weighted average of the resulting state ensemble gives the prior state estimate $\mathbf{x^{-}}(k)$ and the weighted average of the observed ensemble is the model-predicted observation $\mathbf{y}^{-}(k)$.  Covariance matrices $\mathbf{P^{-}}(k)$ and $\mathbf{P^y}(k)$ of the resulting state and observed ensemble, and the cross-covariance matrix $\mathbf{P^{xy}}(k)$ between the state and observed ensembles,  are formed and the equations
\begin{eqnarray} \label{e3}
\mathbf{K}(k) &=& \mathbf{P^{xy}}(k)\left(\mathbf{P^{y}}(k)\right)^{-1}\nonumber\\
\mathbf{P^{+}}(k) &=& \mathbf{P^{-}}(k)-\mathbf{P}^{xy}(k)\left(\mathbf{P}^{y}(k)\right)^{-1}\mathbf{P}^{yx}(k)\nonumber\\
\mathbf{x}^{+}(k) &=& \mathbf{x}^{-}(k)+\mathbf{K}(k)\left(\mathbf{y}(k)-\mathbf{y}^{-}(k) \right).
\end{eqnarray}
are used to update the state and covariance estimates with the observation $\mathbf{y}(k)$.

The UKF algorithm described above can be extended to include the \emph{joint estimation} problem allowing for parameter estimation. In this framework, the parameters $\mathbf{p}$ are considered as auxiliary state variables with trivial dynamics, namely $\mathbf{p}_{k+1} = \mathbf{p}_k$. An augmented $n+l$ dimensional state vector can then be formed consisting of the original $n$ state variables and $l$ model parameters allowing for simultaneous state and parameter estimation \cite{voss,sitz2002}.

To implement parametric prediction, the UKF is used to process the training data and obtain an estimate of $\mathbf{p}$, as well as the state at the end of the training set, $\mathbf{x}(T)$. The parameter values are fixed and Eq. \ref{e1} is forward solved from $t(T)$ to generate predictions of the system dynamics over the prediction interval $\left[t(T+1),t(T+T_F)\right]$. Namely, predictions $\textbf{x}(T+1),\textbf{x}(T+2),\hdots,\textbf{x}(T+T_F)$ are calculated.

\subsection*{Takens' Method for Nonparametric Prediction}
Instead of using the mechanistic model described by Eq. \ref{e1}, the system can be represented nonparametrically. Without loss of generality consider the observed variable $x_{j}$. Using Takens' theorem \cite{takens,SYC}, the $d+1$ dimensional delay coordinate vector $x_j^d(T) = \left[x_j(T),x_j(T-\tau),x_j(T-2\tau),\hdots x_j(T-d\tau)\right]$ is formed which represents the state of the system at time $t(T)$. Here $d$ is the number of delays and $\tau$ is the time-delay. 

The goal of nonparametric prediction is to utilize the training data in the interval $\left[t(0),t(T) \right]$ to build local models for predicting the dynamics over the interval $\left[t\left(T+1\right),t\left(T+T_F\right)\right]$. Here, the method of {\it direct prediction} is chosen. Prior to implementation of the direct prediction, a library of delay vectors is formed from the training data of $x_j$. 

Direct prediction begins by finding the $\kappa$ nearest neighbors, as a function of Euclidean distance, to the current delay-coordinate vector $x_j^d(T)$ within the library of delay vectors. Neighboring delay vectors
\begin{eqnarray*}
x_j^d(T') &=& \left[x_j(T'),x_j(T'-\tau),x_j(T'-2\tau),\hdots x_j(T'-d\tau)\right]\\
x_j^d(T'') &=& \left[x_j(T''),x_j(T''-\tau),x_j(T''-2\tau),\hdots x_j(T''-d\tau)\right]\\
&\vdots&\\
x_j^d(T^\kappa) &=& \left[x_j(T^\kappa),x_j(T^\kappa-\tau),x_j(T^\kappa-2\tau),\hdots x_j(T^\kappa-d\tau)\right]
\end{eqnarray*}
are found within the training data and the known $x_j(T'+i), x_j(T''+i), \ldots, x_j(T^\kappa+i)$ points are used in a local model to predict the unknown value $x_j(T+i)$ where $i = 1,2,\hdots,T_F$. In this article, a locally constant model is chosen
\begin{eqnarray}
\label{localconstant}
x_j(T+i) \approx w_j'x_j(T'+i) + w_j''x_j(T''+i) + \hdots + w_j^{\kappa}x_j(T^\kappa+i)
\end{eqnarray}
where $w_j',w_j'',\hdots,w_j^\kappa$ are the weights for the $j^{th}$ state that determine the contribution of each neighbor in building the prediction. In its simplest form, Eq. \ref{localconstant} is an average of the nearest neighbors where $w_j' = w_j'' = \hdots = w_j^\kappa = \frac{1}{\kappa}$. More sophisticated weighting schemes can be chosen, for example assigning the weights based on the Euclidean distance from each neighbor to the current delay vector \cite{perretti,perretti2,hamilton2}. Selection of values for $d$, $\tau$ and $\kappa$ is necessary for implementation of the direct prediction algorithm. These values were optimized, within each example, to give the lowest prediction error (results not shown).

The accuracy of the predicted $x_j(T+i)$ is subject to several factors.  The presence of noise in the training data plays a substantial role in decreasing prediction accuracy. However, recent advancements in nonparametric analysis have addressed the problem of filtering time series without use of a mechanistic model. In \cite{hamilton2016}, a nonparametric filter was developed which merged Kalman filtering theory and Takens' method. The resulting Kalman-Takens filter was demonstrated to be able to reduce significant amounts of noise in data. Application of the method was extended in \cite{hamiltonEPJ} to the case of filtering stochastic variables without a model. In the results presented below, the training data used for nonparametric prediction are filtered first using the method of \cite{hamilton2016,hamiltonEPJ}.

\subsection*{Hybrid Modeling and Prediction: Merging Parametric and Nonparametric Methods}

As an alternative to the parametric and nonparametric methods described above, we propose a hybrid approach which blends the two methods together. In this framework, we assume that a full mechanistic model as described by Eq. \ref{e1} is available. However, rather than using the full model, a subset of the mechanistic equations are used and the remainder of the variables are represented nonparametrically using delay-coordinates.

In formulating this method it is convenient to first think of Eq. \ref{e1} without vector notation
\begin{eqnarray}
\label{e2}
x_1(k+1) &=& f_1\left(t(k),x_1(k),x_2(k),\hdots,x_n(k),p_1,p_2,\hdots,p_l\right)\nonumber\\
x_2(k+1) &=& f_2\left(t(k),x_1(k),x_2(k),\hdots,x_n(k),p_1,p_2,\hdots,p_l\right)\nonumber\\
&\vdots&\\
x_n(k+1) &=& f_n\left(t(k),x_1(k),x_2(k),\hdots,x_n(k),p_1,p_2,\hdots,p_l\right)\nonumber
\end{eqnarray}
Now assume only the first $n-1$ equations of Eq. \ref{e2} are used to model state variables $x_1,x_2,\ldots,x_{n-1}$, while $x_{n}$ is described nonparametrically
\begin{eqnarray}
\label{hybrid}
x_1(k+1) &=& f_1\left(t(k),x_1(k),x_2(k),\hdots,x_{n-1}(k),x_n(k),p_1,p_2,\hdots,p_l\right)\nonumber\\
x_2(k+1) &=& f_2\left(t(k),x_1(k),x_2(k),\hdots,x_{n-1}(k),x_n(k),p_1,p_2,\hdots,p_l\right)\nonumber\\
&\vdots&\\
x_{n-1}(k+1) &=& f_{n-1}\left(t(k),x_1(k),x_2(k),\hdots,x_{n-1}(k),x_n(k),p_1,p_2,\hdots,p_l\right) \nonumber \\
x_{n}(k+1) &\approx& w_{n}'\tilde{x}_n(T'+k+1) + w_{n}''\tilde{x}_n(T''+k+1) + \hdots + w_{n}^{\kappa}\tilde{x}_n(T^\kappa+k+1) \nonumber
\end{eqnarray}
We refer to Eq. \ref{hybrid} as the {\it hybrid model}. Note, in Eq. \ref{e2} only $x_n$ is assumed to be advanced nonparametrically. This is done purely for ease of presentation and the hybrid model can instead contain several variables whose equations are replaced by nonparametric advancement.

The hybrid model has several distinguishing features. Notice, in this framework nonparametrically advanced dynamics are incorporated into mechanistic equations, essentially merging the two lines of mathematical thought. Furthermore, equations for state variables within Eq. \ref{e2} can be replaced only if there are observations which map directly to them, otherwise their dynamics can not be nonparametrically advanced. Finally, the process of replacing equations in the hybrid method will generally result in a reduction in the number of unknown model parameters to be estimated.

In this hybrid scheme, obtaining an estimate of the unknown parameters in the $n-1$ mechanistic equations and an estimate of $\textbf{x}(T)$ requires a combination of the nonparametric analysis developed in \cite{hamilton2016} and traditional parametric methodology. The state variable $x_{n}$, which is not defined by a mechanistic equation in Eq. \ref{e2}, is represented by delay coordinates within the UKF. Therefore at step $k$ we have the hybrid state
\begin{eqnarray*}
\mathbf{x}^{\mathsmaller H}(k) = \left[x_1(k),x_2(k) ,\ldots, x_{n-1}(k),x_n(k), x_n(k-\tau) , x_n(k-2\tau),\ldots , x_n(k-d\tau)\right]^{\mathsmaller T}
\end{eqnarray*}
The UKF as described above is implemented with this hybrid state $\mathbf{x}^{\mathsmaller H}(k)$ and the model described by Eq. \ref{hybrid}. Notice that in the case of the hybrid model when we have to advance the state dynamics and form the prior estimate in the UKF, the advancement is done parametrically for the first $n-1$ states and nonparametrically for the $n^{th}$ state. Similarly to before, we can augment $\mathbf{x}^{\mathsmaller H}$ with the unknown parameters in the $n-1$ mechanistic equations allowing for simultaneous parameter estimation.

Once the training data are processed and an estimate of $\mathbf{x}^{\mathsmaller H}(T)$ and the parameters are obtained, the hybrid model in Eq. \ref{hybrid} is implemented to generate predictions $\mathbf{x}^{\mathsmaller H}(T+1), \mathbf{x}^{\mathsmaller H}(T+2),\ldots, \mathbf{x}^{\mathsmaller H}(T+T_F)$.

\begin{figure}[!ht]
\begin{center}
\includegraphics[width = \columnwidth]{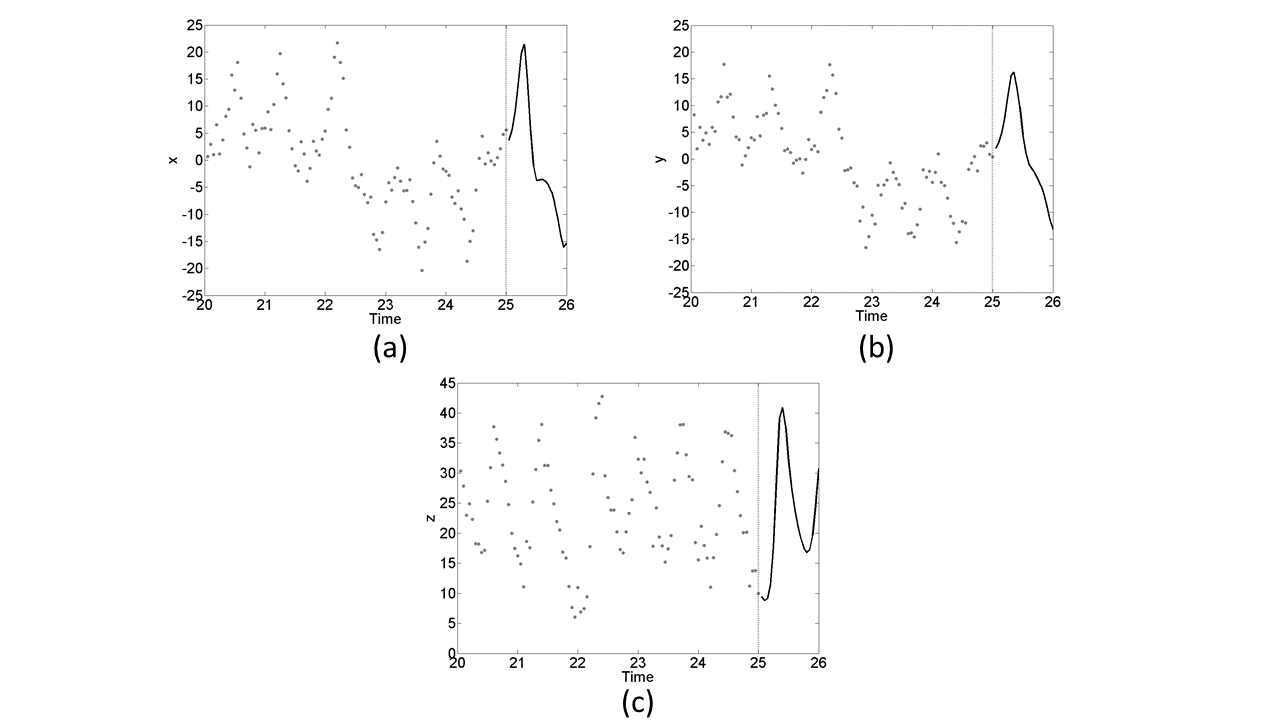}
\end{center}
\caption{{\bf Example of Lorenz-63 realization.} 500 samples, or 25 units of time, of noisy training data (grey circles) are available for (a) $x$, (b) $y$ and (c) $z$. Note, we have only shown the last 5 units of time worth of training data for visualization purposes. From the end of the training data (indicated by dashed black line), we want to accurately predict the system dynamics within the next unit of time (solid black line).}
\label{figure1}
\end{figure}
\section*{Results}
We demonstrate the utility of the hybrid methodology, with comparison to standard parametric and nonparametric modeling and prediction, in the following example systems. When conducting this analysis, two types of error are considered. The first, error in the observations, manifests itself as noise in the training data which all three methods will have to confront. The second type, error in the parameters, takes the form of an uncertainty in the initial parameter values used by the UKF for parameter estimation. Only the parametric and hybrid methods will have to deal with this parameter error. Throughout, we will refer to a percentage uncertainty which corresponds to the standard deviation of the distribution from which the initial parameter value is drawn relative to the mean. For example, if the true value for a parameter $p_1$ is 12 and we have 50\% uncertainty in this value, then the initial parameter value used for estimating $p_1$ will be drawn from the distribution $N(12,(0.5*12)^2)$.

To quantify prediction accuracy, the normalized root-mean-square-error, or SRMSE, is calculated for each prediction method as a function of forecast horizon. Normalization is done with respect to the standard deviation of the variable as calculated from the training data. In using the SRMSE metric, the goal is to be more accurate than if the prediction was simply the mean of the training data (corresponding to SRMSE = 1). Thus a prediction is better than a naive prediction when SRMSE $<$ 1, though for chaotic systems prediction accuracy will eventually converge to this error level since only short-term prediction is possible.

\subsection*{Prediction in the Lorenz-63 System}

As a demonstrative example, consider the Lorenz-63 system \cite{lorenz63}
\begin{eqnarray} \label{lorenz}
\dot{x} &=& \sigma(y-x)\nonumber\\
\dot{y} &=& x(\rho-z)-y\\
\dot{z} &=& xy-\beta z \nonumber
\end{eqnarray}
where $\sigma = 10$, $\rho = 28$, $\beta = 8/3$. Data are generated from this system using a fourth-order Adams-Moulton method with sample rate $h = 0.05$. We assume that 500 training data points of the $x$, $y$ and $z$ variables are available, or 25 units of time. The Lorenz-63 system oscillates approximately once every unit of time, meaning the training set consists of about 25 oscillations. The goal is to accurately predict the dynamics of $x$, $y$ and $z$ one time unit after the end of the training set. However,  the observations of each variable are corrupted by Gaussian observational noise with mean zero and variance equal to 4. Additionally the true value of parameters $\sigma$, $\rho$ and $\beta$ are unknown. Fig. \ref{figure1} shows an example simulation of this system.

The parametric method utilizes Eq. \ref{lorenz} to estimate the model state and parameters, and to predict the $x$, $y$ and $z$ dynamics. For the nonparametric method, delay coordinates of the variables are formed with $d = 9$ and $\tau = 1$. The local constant model for prediction is built using $\kappa = 20$ nearest neighbors. For the hybrid method, the mechanistic equation governing the dynamics of $y$ are replaced nonparametrically resulting in the reduced Lorenz-63 model
\begin{eqnarray*}\label{lorenzhybrid}
\dot{x} &=& \sigma(y-x)\\
\dot{z} &=& xy-\beta z
\end{eqnarray*}
Note, the hybrid model does not require estimation of the $\rho$ parameter since the mechanistic equation for $y$ is removed.

\begin{figure}[!h]
\begin{center}
\includegraphics[width = \columnwidth]{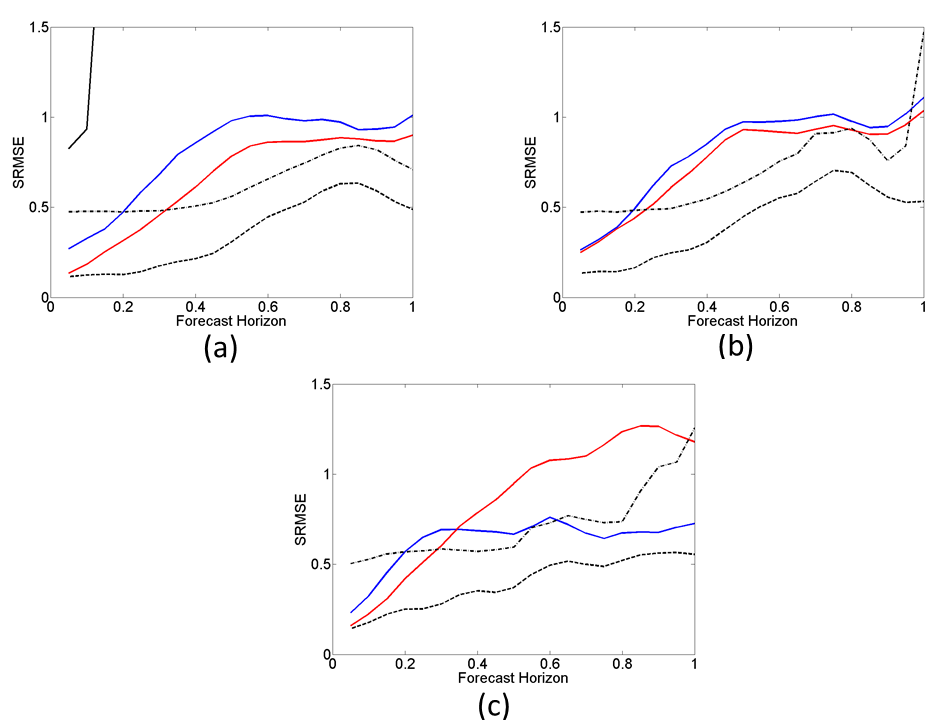}
\end{center}
\caption{{\bf Comparison of the prediction methods in the Lorenz-63 system.}
Results of predicting the Lorenz-63 (a) $x$, (b) $y$ and (c) $z$ variables averaged over 500 realizations. Training data consists of 500 data points generated from Eq. \ref{lorenz} with $\sigma = 10$, $\rho = 28$ and $\beta = 8/3$ with sample rate $h = 0.05$. Data are corrupted by Gaussian observational noise with mean 0 and variance of 4. Parametric (black), nonparametric (blue) and hybrid (red) prediction SRMSE plotted as a function of forecast horizon with initial parameter uncertainty of 80\% (solid line), 50\% (dashed-dotted line) and 20\% (dashed line). Hybrid prediction with 80\% uncertainty, which utilizes mechanistic equations in describing $x$ and $z$ but nonparametrically represents $y$, offers an improvement in short-term prediction accuracy over standalone nonparametric prediction. Parametric prediction with 80\% uncertainty performs poorly for predicting all three variables and in the case of (b) and (c) is not seen due to the scale of the error. As the uncertainty shrinks, performance of the parametric method improves. However, only at a small uncertainty level does the parametric method outperform the short-term improvement in prediction afforded by the hybrid method.}
\label{figure2}
\end{figure}

Fig. \ref{figure2} shows a comparison of parametric (black), nonparametric (blue) and hybrid (red) prediction error as a function of forecast horizon. SRMSE results averaged over 500 system realizations. Various parameter uncertainty levels are shown: 80\% uncertainty (solid lines), 50\% uncertainty (dashed-dotted lines) and 20\% uncertainty (dashed line). The hybrid method with 80\% uncertainty offers improved short-term prediction of the Lorenz-63 $x$ (Fig. \ref{figure2}a) and $z$ (Fig. \ref{figure2}c) variables over standalone nonparametric prediction as well as parametric prediction with 80\% uncertainty. Hybrid and nonparametric prediction of $y$ (Fig. \ref{figure2}b) are comparable, which is to be expected since the hybrid approach is using nonparametric advancement of $y$ in its formulation. Note that parametric prediction at this uncertainty level does very poorly and in the cases of $y$ and $z$ its result is not shown due to the scale of the error. As the uncertainty decreases for parametric prediction, its performance improves. However, hybrid prediction with 80\% uncertainty still outperforms parametric prediction with 50\% uncertainty in the short-term. At a small uncertainty level, parametric prediction outperforms both hybrid and nonparametric methods which is to be expected since it has access to the true model equations and starts out with close to optimal parameter values.

The success of the hybrid method at higher uncertainty levels can be traced to more accurate estimates of the model parameters in the mechanistic equations that it uses. Table \ref{table_L63} shows the resulting hybrid and parametric estimation of the Lorenz-63 parameters. The hybrid method with 80\% uncertainty is able to construct accurate estimates of both $\sigma$ and $\beta$, with a mean close to the true value and a small standard deviation of the estimates. The parametric method with 80\% and 50\% uncertainty is unable to obtain reliable estimates, exemplified by the large standard deviation of the estimates. Only when the parametric method has a relatively small uncertainty of 20\% is it able to accurately estimate the system parameters.

\begin{table}[ht]
\begin{center}
\begin{tabular}{| c | c | c | c | c |}
\hline
\multicolumn{4}{|c|}{Lorenz-63 Parameter Estimation Results} \\ \hline
True Parameter & Method&  Mean&Standard Deviation \\ \hline
\multirow{4}{*}{$\sigma = 10$} & Hybrid (80\% Uncertainty) & 9.77 & 0.75 \\
 & Parametric (80\% Uncertainty) & 8.03 & 4.81 \\
& Parametric (50\% Uncertainty) & 9.84 & 3.06  \\ 
& Parametric (20\% Uncertainty) & 10.06 & 0.95 \\
\cline{1-4}
\multirow{4}{*}{$\rho = 28$} & Hybrid (80\% Uncertainty) & NA & NA\\
 & Parametric (80\% Uncertainty) & 24.55 & 14.07 \\
& Parametric (50\% Uncertainty) & 25.63 & 6.37  \\ 
& Parametric (20\% Uncertainty) & 27.89 & 0.83 \\
\cline{1-4}
\multirow{3}{*}{$\beta = 2.67$} &  Hybrid (80\% Uncertainty) & 2.58 & 0.11\\
 & Parametric (80\% Uncertainty) & 1.61 & 1.34 \\
& Parametric (50\% Uncertainty) & 2.20 & 0.98  \\ 
& Parametric (20\% Uncertainty) & 2.63 & 0.19 \\
\cline{1-4}
\end{tabular}
\caption{\textbf{Summary of Lorenz-63 parameter estimation results}. Mean and standard deviation calculated over 500 realizations. The hybrid method, which only needs to estimate $\sigma$ and $\beta$, is robust to a large initial parameter uncertainty. The parametric method on the other hand is unable to obtain reliable estimates of the Lorenz-63 parameters unless the uncertainty is small enough.}
\label{table_L63}
\end{center}
\end{table}

\subsection*{Predicting Neuronal Network Dynamics}
We now consider the difficult high dimensional estimation and prediction problem posed by neuronal network studies. If we are only interested in predicting a portion of the network, then we can use the proposed hybrid method to refine our estimation and prediction while simultaneously reducing estimation complexity. As an example
of this potential network application we consider the prediction of spiking dynamics in a network of $M$ Hindmarsh-Rose neurons \cite{hindmarsh}
\begin{eqnarray}\label{hindmarsh}
\dot{x}_i &=& y_i-a_ix_i^3+b_ix_i^2-z_i+1.2+\sum_{i\neq m}^M \frac{\beta_{im}}{1+9e^{-10x_m}}x_m \nonumber \\
\dot{y}_i &=& 1-c_ix_i^2 \\
\dot{z}_i &=& 5\times 10^{-5}\left[4\left(x_i-\left(-\frac{8}{5}\right)\right)-z_i \right] \nonumber
\end{eqnarray} 
where $i = 1,2,\hdots,M$. $x_{i}$ corresponds to the spiking potential while $y_i$ and $z_i$ describe the fast and slow-scale dynamics, respectively, of neuron $i$. Each individual neuron in the network has parameters $a_i =1, b_i = 3$ and $c_i = 5$ which are assumed to be unknown. $\beta_{im}$ represents the connectivity coefficient from neuron $i$ to neuron $m$. For a network of size $M$, we have $M^2-M$ possible connection parameters since neuron self connections are not allowed (i.e. $\beta_{ii} = 0$). These connection parameters are also assumed to be unknown.

\begin{figure}[ht]
\begin{center}
\includegraphics[width = \columnwidth]{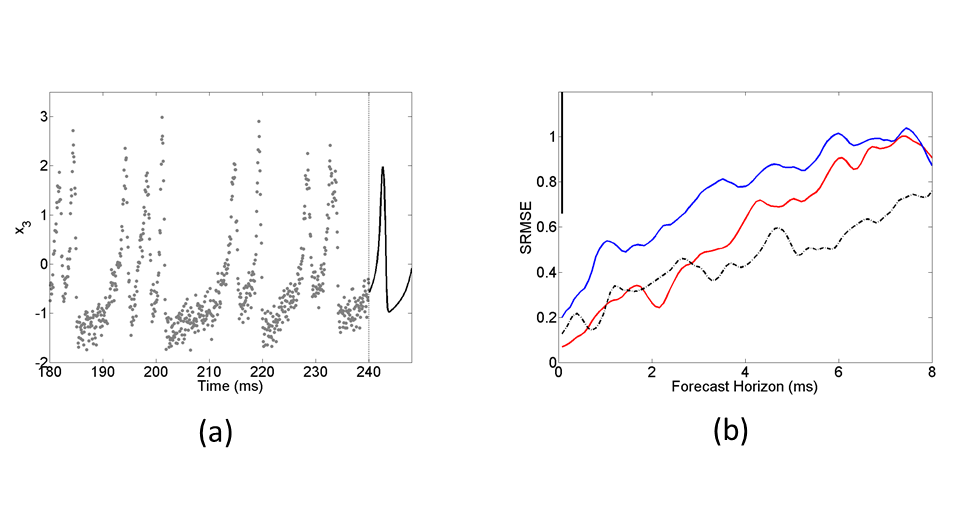}
\end{center}
\caption{{\bf Predicting neuron potential $x_3$ in random 3-neuron Hindmarsh-Rose networks.} (a) 3000 samples (or 240 ms) of training data (grey circles) are available from each neuron in the network. From the end of the training data (indicated by dashed black line), we want to accurately predict the next 8 ms of $x_3$ (solid black line). (b) Average SRMSE when predicting 8 ms of $x_3$. Results averaged over 200 randomly generated 3-neuron Hindmarsh-Rose network realizations. Prediction accuracy when using the full parametric model (black), nonparametric (blue) and the hybrid method (red) shown. 80\% uncertainty (solid line) and 50\% uncertainty (dashed-dotted line) levels shown. Once again, the hybrid method with 80\% uncertainty offers improved accuracy in predicting $x_3$ over the nonparametric and parametric with 80\% uncertainty methods. Prediction accuracy between the hybrid and parametric with 50\% uncertainty is comparable.
}
\label{figure3}
\end{figure}

For this example we examine networks of size $M = 3$ with 5 random connections. Data from these networks are generated using a fourth-order Adams-Moulton method with sample rate $h = 0.08$ ms.  We assume that the training data consists of 3000 observations, or 240 ms, of the $x_1, x_2,x_3$ variables each of which are corrupted by Gaussian noise with mean 0 and variance of 0.2. Under the stated parameter regime, the neurons in the network spike approximately every 6 ms, meaning our training set has on average around 40 spikes per neuron. In this example, we restrict our focus to predicting 8 ms of the $x_3$ variable (though a similar analysis follows for the prediction of $x_1$ and $x_2$). Fig. \ref{figure3}a shows a representative realization of this problem. Given our interest in $x_3$, the hybrid method only assumes a mechanistic equation for neuron 3
\begin{eqnarray*}
\dot{x}_3 &=& y_3-a_3x_3^3+b_3x_3^2-z_3+1.2+\sum_{3 \neq m}^M \frac{\beta_{3m}}{1+9e^{-10x_m}}x_m \\
\dot{y}_3 &=& 1-c_3x_3^2\\
\dot{z}_3 &=& 5\times 10^{-5}\left[4\left(x_3-\left(-\frac{8}{5}\right)\right)-z_3 \right]
\end{eqnarray*}
and nonparametrically represents neuron 1 and neuron 2.

Fig. \ref{figure3}b shows the resulting accuracy in predicting $x_3$ when using parametric (black), nonparametric (blue) and hybrid (red) methods with 80\% (solid line) and 50\% (dashed-dotted line) uncertainty in parameter values. The parametric approach uses the full mechanistic model described by Eq. \ref{hindmarsh} for modeling and prediction, requiring estimation of the $x,y$ and $z$ state variables and parameters $a,b$ and $c$ for each neuron, as well as the full connectivity matrix. Notice that once again with 80\% uncertainty, the scale of error for the parametric method is much larger compared to the other methods. Only with 50\% uncertainty is the parametric method able to provide reliable predictions of $x_3$. Note that unlike in the Lorenz-63 example, we do not consider the parametric method with 20\% uncertainty since reasonable parameter estimates and predictions are obtained with 50\% uncertainty. The nonparametric method ($\tau = 1$, $d = 9$) uses $\kappa = 10$ neighbors for building the local model for prediction. Again we observe that the hybrid method, even with a large parameter uncertainty of 80\%, provides accurate predictions of $x_3$ compared to the other methods. Table \ref{table_HR} shows the robustness of the hybrid method in estimating the individual parameters for neuron 3.

\begin{table}[ht]
\begin{center}
\begin{tabular}{| c | c | c | c | c |}
\hline
\multicolumn{4}{|c|}{Neuron 3 Parameter Estimation Results} \\ \hline
True Parameter & Method&  Mean&Standard Deviation \\ \hline
\multirow{4}{*}{$a_3 = 1$} & Hybrid (80\% Uncertainty) & 0.98 & 0.04 \\
 & Parametric (80\% Uncertainty) & 1.07 & 0.51 \\
& Parametric (50\% Uncertainty) & 0.98 & 0.15  \\ 
\cline{1-4}
\multirow{4}{*}{$b_3 = 3$} & Hybrid (80\% Uncertainty) & 2.96 & 0.10\\
 & Parametric (80\% Uncertainty) & 2.92 & 0.88 \\
& Parametric (50\% Uncertainty) & 2.92 & 0.26  \\ 
\cline{1-4}
\multirow{3}{*}{$c_3 =5$} &  Hybrid (80\% Uncertainty) & 4.93 & 0.16\\
 & Parametric (80\% Uncertainty) & 4.66 & 1.04 \\
& Parametric (50\% Uncertainty) & 4.83 & 0.43  \\ 
\cline{1-4}
\end{tabular}
\caption{\textbf{Summary of neuron 3 parameter estimation results}. Mean and standard deviation calculated over 200 realizations. The hybrid method once again is robust to a large initial parameter uncertainty. The parametric method on the other hand is unable to obtain reliable estimates of the neuron parameters with large uncertainty.}
\label{table_HR}
\end{center}
\end{table}

\subsection*{Predicting Flour Beetle Population Dynamics}

We now investigate the prediction problem in a well-known data set from an ecological study involving the cannibalistic red flour beetle \emph{Tribolium castaneum}. In \cite{constantino}, the authors present experimentally collected data and a mechanistic model describing the life cycle dynamics of \emph{T. castaneum}. Their discrete time model describing the progression of the beetle through the larvae, pupae, and adult stages is given by 
\begin{eqnarray} \label{beetle}
L(t+1) &=& bA(t) e^{-c_{el}L(t) - c_{ea}A(t)} \nonumber\\
P(t+1) &=& L(t)(1-\mu_l)\\
A(t+1) &=& P(t) e^{-c_{pa}A(t)}+A(t)(1-\mu_a)\nonumber
\end{eqnarray}
where $L$, $P$ and $A$ correspond to larvae, pupae and adult populations, respectively. The essential interactions described by this model are (i) flour beetles become reproductive only in the adult stage, (ii) adults produce new larvae, (iii) adults and larvae can both cannibalize larvae, and (iv) adults cannibalize pupae. We note that since Eq. \ref{beetle} only approximates the life cycle dynamics of the beetle, there is a degree of model error in the proposed system, unlike the previous examples.

\begin{figure}[!h]
\begin{center}
\includegraphics[width = \columnwidth]{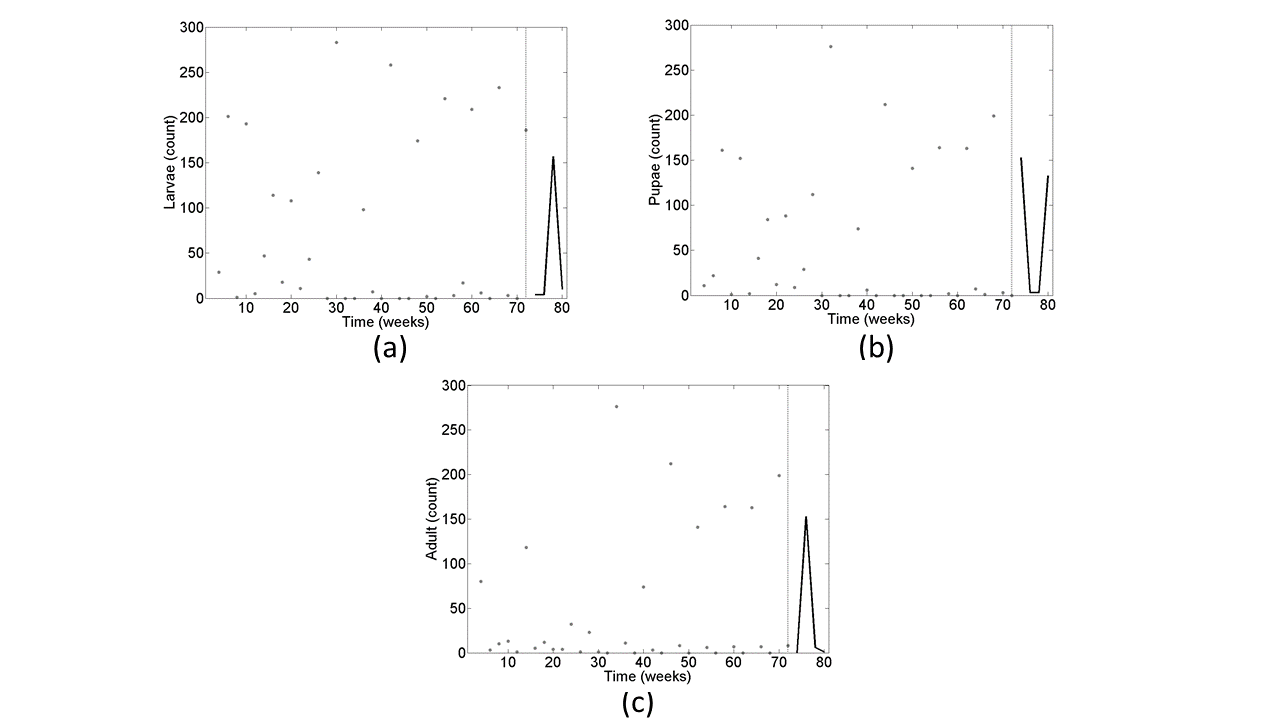}
\end{center}
\caption{{\bf Example data set from \emph{T. castaneum} experiment presented in \cite{constantino}.}
36 observations, or 72 weeks, of training data (grey circles) are available for (a) larvae, (b) pupae and (c) adult population levels. From the end of the training data (indicated by dashed black line), we want to accurately predict the next 8 weeks of population dynamics (solid black line).}
\label{figure4}
\end{figure}

The authors of \cite{constantino} experimentally set the adult mortality rate ($\mu_a$) to $0.96$ and the recruitment rate ($c_{pa}$) from pupae to adult to seven different values ($0$, $0.05$, $0.10$, $0.25$, $0.35$, $0.50$, $1.0$). Experiments at each recruitment rate value were replicated three times resulting in 21 different datasets. Each dataset consists of total numbers of larvae, pupae, and adults measured bi-weekly over 82 weeks resulting in 41 measurements for each life stage. These data were fit to Eq. \ref{beetle} in \cite{constantino} and parameter estimates $b = 6.598$, $c_{el} = 1.209 \times 10^{-2}$, $c_{ea} = 1.155 \times 10^{-2}$ and $\mu_l = 0.2055$ were obtained. We treat these parameter values as ground truth when considering the different parameter uncertainty levels for fitting the data to the model.

\begin{figure}[ht]
\begin{center}
\includegraphics[width = \columnwidth]{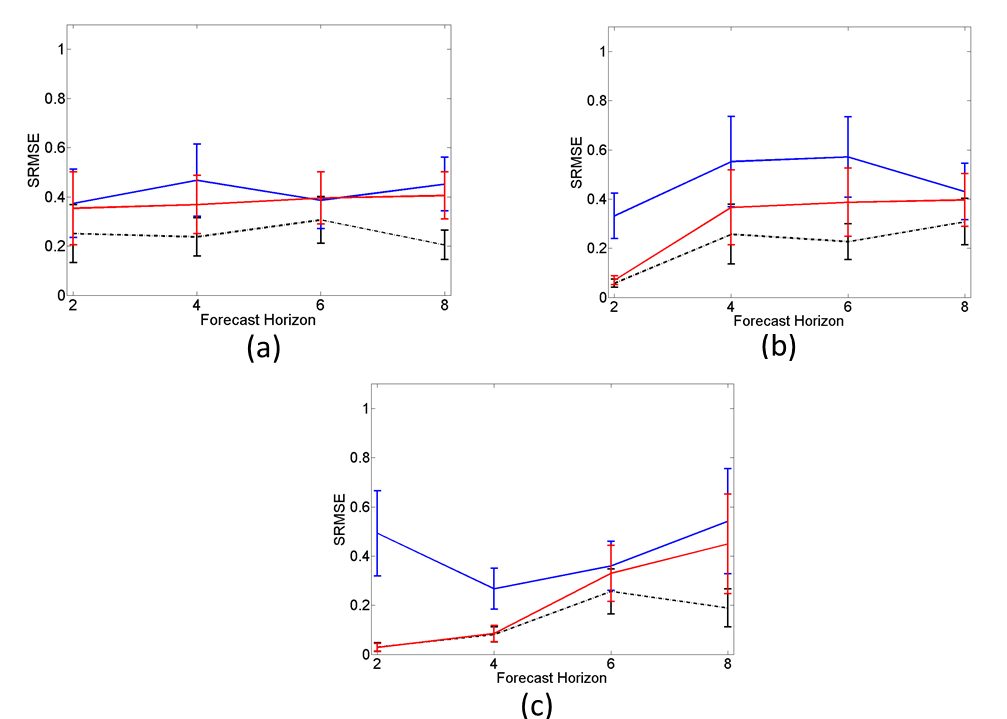}
\end{center}
\caption{{\bf Results for predicting population levels of \emph{T. castaneum}.} Average SRMSE over 21 experimental datasets when using parametric (black curve), nonparametric (blue curve) and hybrid (red curve) methods for predicting (a) larvae, (b) pupae and (c) adult population levels with uncertainty of 80\% (solid line) and 50\% (dashed-dotted line). Error bars correspond to standard error over the 21 datasets. Hybrid prediction with 80\% uncertainty offers improved prediction over both nonparametric and parametric with 80\% uncertainty (not visible due to scale of error), and comparable performance to parametric with 50\% uncertainty.}
\label{figure5}
\end{figure}

In our analysis of this system, we treat the first 37 measurements (or 74 weeks) within an experiment as training data and use the remaining 4 time points (or 8 weeks) for forecast evaluation. Fig. \ref{figure4} shows an example of this setup for a representative dataset. Fig. \ref{figure5} shows the results of predicting the larvae (Fig. \ref{figure5}a), pupae (Fig. \ref{figure5}b) and adult (Fig. \ref{figure5}c) populations using parametric (black), nonparametric (blue) and hybrid prediction methods with 80\% (solid line) and 50\% (dashed-dotted line) parameter uncertainty levels. Error bars correspond to the standard error over the 21 datasets. The parametric method uses the full mechanistic model described in Eq. \ref{beetle} to estimate the population state and parameters $b, c_{el}, c_{ea}$ and $\mu_l$ before prediction. We note in Fig. \ref{figure5} that the parametric method with 80\% uncertainty is not shown due to the scale of the error, and is significantly outperformed by the nonparametric prediction ($\tau = 1, d = 2, \kappa = 5$). For the hybrid method, we only consider the mechanistic equations for pupae and adult population dynamics
\begin{eqnarray*}
P(t+1) &=& L(t)(1-\mu_1)\\
A(t+1) &=& P(t) e^{-c_{pa}A(t)}+A(t)(1-\mu_a)
\end{eqnarray*}
and nonparametrically represent larvae. Hybrid prediction with 80\% uncertainty outperforms both nonparametric and parametric with 80\% uncertainty for pupae and adult population levels, and is comparable to parametric with 50\% uncertainty.

\section*{Conclusion}

By blending characteristics of parametric and nonparametric methodologies, the proposed hybrid method for modeling and prediction offers several advantages over standalone methods. From the perspective of model fitting and the required parameter estimation that arises in this process, we have shown that the hybrid approach allows for a more robust estimation of model parameters. Particularly for situations where there is a large uncertainty in the true parameter values, the hybrid method is able to construct accurate estimates of model parameters when the standard parametric model fitting fails to do so. At first this may seem counter-intuitive, but in fact it is not that surprising. The replacement of mechanistic equations with their nonparametric representations in effect reduces the dimension of the parameter space that we have to optimize in, resulting in better parameter estimates. As we have demonstrated in the above examples, this refinement in the parameter estimates leads to an improvement in short-term prediction accuracy.

The limitations of the hybrid method are similar to those of parametric and nonparametric methods in that if not enough training data are available then accurate estimation and prediction becomes difficult. However, the demonstrated robustness of the hybrid method to large parameter uncertainty is encouraging, particularly when considering experimental situations where we may not have a good prior estimate of the model parameters. One could consider implementing the hybrid method in an iterative fashion, estimating the parameters of each equation separately, then piecing the model back together for prediction. We can think of this as an {\it iterative hybrid method}, and is the subject of future work.

We view this work as complementary to recent publications on forecasting \cite{perretti,perretti2,hartig}. The authors of \cite{perretti,perretti2} advocate nonparametric methods over parametric methods in general, while a letter \cite{hartig} addressing the work of \cite{perretti} showed that a more sophisticated method for model fitting results in better parameter estimates and therefore model-based predictions which outperform model-free methods. Our results support the view that no one method is uniformly better than the other. As we showed in the above examples, in situations where the model error and uncertainty in initial parameters are relatively small, the parametric approach outperforms other prediction methods. Often in experimental studies though, we are not operating in this ideal situation and instead are working with a model that has substantial error with a large uncertainty in parameters which can lead to inaccurate system inference. In situations such as these, nonparametric methods are particularly useful.

The main appeal of the hybrid method is that we can confront these situations without having to completely abandon the use of the mechanistic equations. This is important since mechanistic models often provide valuable information about the underlying processes governing the system dynamics. While we explored in detail the robustness of the hybrid method to large levels of parameter uncertainty, its usefulness stretches well beyond that. In some instances, we may only have a model for some of the states or portions of the model may have higher error than others. By supplementing these parts with their nonparametric representation, the hybrid method would allow us to only use the parts of the model we are confident in and thus improve our analysis.


\nolinenumbers

\providecommand{\noopsort}[1]{}\providecommand{\singleletter}[1]{#1}%

\end{document}